\documentclass{article}
\usepackage[english]{babel}
\usepackage{amsmath,amssymb}

\newcommand{\nobracket}{}
\newcommand{\tmem}[1]{{\em #1\/}}
\newcommand{\tmmathbf}[1]{\ensuremath{\boldsymbol{#1}}}
\newcommand{\tmop}[1]{\ensuremath{\operatorname{#1}}}
\newcommand{\tmrsup}[1]{\textsuperscript{#1}}
\newcommand{\tmstrong}[1]{\textbf{#1}}

\begin{document}

\begin{center}
  {\tmstrong{PERSISTENCE OF GAPS IN THE SPECTRUM OF CERTAIN}}
\end{center}

\begin{center}
  {\tmstrong{ALMOST PERIODIC OPERATORS}}
\end{center}

\

\begin{center}
  {\tmstrong{Norbert Riedel}}
\end{center}

\

{\tmem{Abstract:}} It is shown that for any irrational rotation number and
any admissible gap labelling number the almost Mathieu operator (also known as
Harper's operator) has a gap in its spectrum with that labelling number. This
answers the strong version of the so-called ``Ten Martini Problem''. When
specialized to the particular case where the coupling constant is equal to
one, it follows that the ``Hofstadter butterfly'' has for any quantum Hall
conductance the exact number of components prescribed by the recursive scheme
to build this fractal structure.

.

{\tmstrong{Introduction}}

The present work is concerned with the spectral properties of the simplest
kind of discrete Schroedinger type operators with an almost periodic
potential. These operators form a self-dual class with respect to the Fourier
transform, with exactly one operator being invariant. More specifically, a
precise one-to-one relationship between gaps and admissible gap labelling
numbers will be established. While relegating the technical formulation of the
problem and the outline of its proof to the first section, we are going to
dwell for the rest of the introduction on the significance of this result in
the special case of the self-dual operator for physics and mathematics.

Recent progress that has occurred in solid state physics through the
development of new experiments which led to improved measurements, has made it
possible to obtain improved physical evidence for the presence of the
butterfly fractal spectrum. In [GSUKNKS] portions of the ``Hofstadter
butterfly'' have been observed in lateral superlattices patterned on
GaAs/AlGaAs heterostructure by exploiting the quantum Hall conductance as a
diagnostic tool. One way to view the significance of the present paper is the
recognition, that this diagnostic approach to detecting the fractal stucture
experimentally is on solid theoretical ground, at least as long as the
``Hofstadter butterfly'' is accepted as a paradigm for the quantum Hall
effect: Knowing the specific value of the quantum Hall conductance, it is
possible to allocate the corresponding gap components in the fractal structure
in a prescribed way, without running the risk of missing one. This follows
from the considerations in [OA], where it was shown that, assuming that the
strong version of the ``Ten Martini problem'' holds, the components associated
with the same quantum Hall conductance can be counted by means of a specific
combinatorial formula. In another experiment cold neutral atoms in an optical
lattice were used to exhibit salient features of the butterfly fractal [JZ].
The employment of lasers allows the simulation of a magnetic flux through the
lattice. The success of this approach for specific rotation numbers depends on
the visibility of the particle density, which in turn depends on its
periodicity. As shown by a specific experiment in [JZ] (FIG 4b), the
visibility decreases for the irrational rotation number $\alpha = \frac{1}{2
\pi}$. Again, as the present work will show, it is reassuring to know, that
the butterfly fractal does not allow for unexpected transitions at irrational
rotation numbers to occur, which could not be picked up through a suitably
designed experiment.

We turn now to the mathematical significance of the present work, which at
this point is somewhat speculative in nature. The set of quantum Hall numbers
is a cyclic subgroup of the additive group of real numbers. If $\tmmathbf{g}$
\ is the posiitve generator of this group, then each quantum Hall number
$\tmmathbf{q}$ can be written as $\tmmathbf{q}= k\tmmathbf{g}$, for a suitable
integer {\tmem{k}}. As was shown in [OA] under the assumption that the strong
version of the ``Ten Martini conjecture'' holds, the number of components in
the butterfly fractal with a common positive quantum Hall number
$\tmmathbf{q}$ is equal to
\[ \Phi (2 k), \tmop{where} \Phi (n) = \sum_{j = 1}^n \varphi (j), \]
and $\varphi$ is Euler's totient function. On the other hand, an observation
by J. Franel from the 1920's asserts that the Riemann hypothesis is true if
and only if the relation
\[ \sum_{j \leqq \Phi (n)} (r^{(n)}_j - \frac{j}{\Phi (n)})^2 = O (n^{- 1 +
   \varepsilon}) \]
holds for any $\varepsilon > 0$. Here $r^{(n)}_j$ is the $j$th Farey fraction
of order n. For a detailed exposition of this subject see Edmund Landau's
lectures on number theory [L], Band 2, Kapitel 13, or [E], section 12.2.
Notice that the Farey fractions, in the order as they appear in this
asymptotic formula, occupy a prominent position in the butterfly fractal. So
it appears that the butterfly fractal holds information about the Riemann
hypothesis, the exact nature of which remains yet to be determined. Inspired
by the color coding of the butterfly fractal according to the Hall
conductance, which was also introduced in [OA], and which has become quite
popular in recent years, it is tempting to take a cue from Marc Kac, who is
not only alleged to have offered ten martinis as a reward for the solution of
the eponymous problem which is the subject of this work, but who also famously
asked ``Can one hear the shape of a drum?'', and pose the question, ``Can one
see the hue of the Riemann hypothesis?''

Finally, we turn to a description of the organization of the paper. In section
1 the problem will be introduced in a form that is conducive to the employment
of tools which are needed to solve it. Based on three propositions, two of
which are known results, while the third one still needs to be established,
the short proof of the major result will be given in that section. In section
2 material from this author's previous work will be assembled in a fashion
that facilitates its usage in the present context, and a number of preparatory
results will be established. In section 3 the proof for the outstanding
proposition will be provided. \ \ \ \ \ \ \ \ \ \ \ \ \ \ \ \ \ \ \ \ \ \ \ \
\ \ \ \ \ \ \ \ \ \ \ \ \ \ \ \ \ \ \ \ \ \ \ \ \ \ \ \

.

{\tmstrong{Section 1}}

The observation that gaps tend to open up readily for small coupling constants
due to basic perturbations of the degenerate case, that is when the coupling
constant is equal to zero, naturally leads to a search for an argument that
allows one to show that those gaps can not close as the coupling constant
increases. The advantage of such an approach goes beyond mere expediency. By
establishing the persistence of gaps as opposed to proving the existence of
gaps for specific parameters, one can reach through to the elusive self-dual
case ``from the outside''. It is the purpose of the present work to develop
such an argument. Earlier attempts can be readily traced in the literature.
For instance the major thrust in [HKS] is along the same line: If one
substitutes the conjecture C3 in that paper, which at this point remains
unproven, by the combination of Proposition 2 and Proposition 3 below, then
one can simply follow the argument provided in ``Remarque 3.2.4'' in [HKS] to
obtain the desired result. In conclusion our proof will require two legs to
stand on: The first leg establishes that gaps with prescribed labels open up
for sufficiently small (but by no means uniform!) coupling constants. This has
been rigorously shown in [HKS]. The second leg guarantees that (open) gaps
can't close, as the coupling constant increases. Our argument to accomplish
this has two components. First, we need to show that if a gap closes, then it
follows that the Lyapunov exponent, considered as a function of two real
parameters, namely the coupling constant and the spectral parameter, has a
local maximum in one of the resolvent sets. Second, we need to prove that the
Lyapunov exponent does not have any critical points. It is the proof of the
second component that will occupy the main body of the paper.

We turn now to outlining the setting for the proof of that particular
component. Let $\alpha$ be an irrational number, let {\tmem{u}} and {\tmem{v}}
be unitary operators satisfying the relation

(1.1)
\[ \text{{\tmem{$u v = e^{2 \pi \alpha i}$}}} v u, \]
and let {\tmem{$A_{\alpha}$}} be the (abstract) $C^{\ast}$-algebra generated
by {\tmem{u}} and {\tmem{v}}. Furthermore, let $\tau$ be the unique tracial
state on $A_{\alpha}$. This is a positive linear functional, standardized by
setting $\tau (\mathbb{e}) = 1$, where $\mathbb{e}$ denotes the unit element
in $A_{\alpha}$, such that $\tau (a b) = \tau (b a$) holds for all elements in
$a, b \epsilon A_{\alpha}$. In this setting we define for any positive
coupling constant $\beta$ the almost Mathieu operator as follows,

(1.2)
\[ h (\beta) = u  + u^{\ast} + \beta (v + v^{\ast}) . \]
As usual, the upper right asterisk denotes the adjoint of an operator. The
integrated density of states can now be identified with the restriction of the
functional $\tau$ to the abelian $C^{\ast}$-algebra generated by $h (\beta)$.
To obtain the integrated density of states proper, all one has to do is to
represent this restricted functional by a probability measure on the spectrum
of $h (\beta)$. In light of the comments made above we can now formulate our
first proposition which establishes the first leg of the argument.

\

{\tmstrong{1.1 Proposition}} ([HKS]) For any number $r \epsilon [0, 1] \cap
\mathbb{Z}+ \alpha \mathbb{Z}$ there exists a positive number $\beta_0$ with
the property that for any $\beta \epsilon (0, \beta_0$], the operator $h
(\beta$) has a gap in its spectrum with the label {\tmem{r}}. More
specifically, there exists a real number $s_{\beta}$ in the resolvent set of
$h (\beta$), such that the spectral projection {\tmem{p}} associated with the
interval $(- \infty, s_{\beta}$] has\, the property $\tau (p) = r$.

\

We turn now to the second leg. First we need to define the Lyapunov exponent
in a way that is compatible with the present settings.
\[ L : \mathbb{R} \times \mathbb{C} \rightarrow \mathbb{R}, L (\beta, z) =
   \tau (\log | h (\beta) - z |) . \]
While the operator $\log | h (\beta) - z | = \log [(h (\beta) - z) (h (\beta)
- \bar{z})]^{\frac{1}{2}}$ is an element of $A_{\alpha}$ for $z \epsilon
\mathbb{C} \backslash S p (h (\beta))$, \ $S p (h (\beta))$ denoting the
spectrum of $h (\beta)$, this is not the case for $z \epsilon S p (h
(\beta))$. However, for any complex number this operator is contained in $L^1
(A, \tau)$, the space of ``integrable operators'' associated with $A_{\alpha}$
and $\tau$. By virtue of the so-called Thouless formula the number $L (\beta,
z)$ is seen to coincide with the usual definition of the Lyapunov exponent for
$h (\beta)$ at {\tmem{z}}. We can now formulate our next statement.

\

{\tmstrong{1.2 Proposition}} The function {\tmem{L}} is jointly continuous is
both variables. Moreover, for $\beta \leq 1$, $L (\beta, z) = 0$ for every $z
\epsilon S p (h (\beta))$.

\

While this follows for rotation numbers satisfying a diophantine condition
from the present author's earlier work in conjunction with the semicontinuity
of the spectrum (for badly approximable numbers see [R4], Proposition 2.12,
and for sufficiently well approximable numbers see [R2], Corollary 2.3), a
proof which does not rely on any diophantine condition is implicit in [BJ].
Indeed, the crucial Proposition 9 towards the end of that paper is valid for
any pair of sufficiently smooth potentials, not just those which differ by the
spectral parameter {\tmem{E}}. The remaining argument then carries through
essentially without change.

The third statement, which will be proved in section 2 and section 3, is as
follows.

\

{\tmstrong{1.3 Proposition}} For every $z \epsilon \mathbb{R} \backslash S p
(h (\beta))$,
\[ (\tau ((h (\beta) - z)^{- 1}), \tau ((h (\beta) - z)^{- 1} v)) \neq (0, 0)
   . \]

The proof of the following theorem is the main objective in the sequel.

\

{\tmstrong{1.4 Theorem}} For every $r \epsilon [0, 1] \cap (\mathbb{Z}+
\alpha \mathbb{Z})$ there exists $s \epsilon \mathbb{R} \backslash S p (h
(\beta))$ such that the spectral projection {\tmem{p}} associated with the
interval $(- \infty, s]$ has the property $\tau (p) = r$.

\

{\tmstrong{Proof}}: \ By \ duality it suffices to consider the case $0 <
\beta \leq 1$ only. Let $r \epsilon [0, 1] \cap (\mathbb{Z}+ \alpha
\mathbb{Z})$, and choose $\beta_0$ as in Proposition 1.1. For some
$\tilde{\beta} \epsilon (0, \beta_0)$ let $s_{\tilde{\beta}}$ be as in
Proposition 1.1. Let $\Omega_r$ be the connected component of $R = \{ (\beta,
t) \epsilon \mathbb{R}^+ \times \mathbb{R}/ t \epsilon \mathbb{R} \backslash S
p (h (\beta)) \}$ containing the point ($\tilde{\beta}$, $s_{\tilde{\beta}}$).
Since $\tilde{\beta}$ can be chosen arbitrarily close to zero, it suffices to
show that $\Omega_r$ contains a point whose first coordinate is equal to 1.
Suppose this were not true. Then $\Omega_r \subset [0, 1] \times [- 4, 4]$,
and therefore $\overline{\Omega }_r$ is compact. By Proposition 1.2 the
function {\tmem{L}} is continuous on $\overline{\Omega }_r$, and it takes the
value zero on the boundary $\partial \overline{\Omega }_r$. It follows that
{\tmem{L}} has a local maximum at some point ($\beta_e, s_e$) in $\Omega_r$.
Since {\tmem{L}} is infinitely differentiable in $\Omega_r$, the gradient of
{\tmem{L}} at ($\beta_e, s_e$) vanishes. Thus
\[ \frac{\partial L}{\partial z} \text{ } (\beta_e, s_e) = \tau ((s_e - h
   (\beta_e))^{- 1}) = 0, \]
\[ \frac{\partial L}{\partial \beta} \text{ } (\beta_e, s_e) = \tau ((h
   (\beta_e) - s_e)^{- 1} (v + v^{\ast})) = 2 \tau ((h (\beta_e) - s_e)^{- 1}
   v) = 0. \]
By Proposition 1.3 this is impossible.

\

{\tmstrong{Remarks:}} 1) While the first of the two partial derivatives
occurring in the proof is obvious, the second one warrants a few words of
explanation. Differentiating the first of the above partial derivatives with
respect to $\beta$, by invoking the formula
\[ \frac{\partial}{\partial \beta} (z - h (\beta))^{- 1} = (z - h (\beta))^{-
   1} (v + v^{\ast}) (z - h (\beta))^{- 1}, \]
and then antidifferentiating the result with respect to {\tmem{z}} yields the
claimed formula plus a function which depends on $\beta$ only, $f (\beta)$
say. In order to show that $f (\beta)$ is actually zero, one would simply like
to let {\tmem{z}} approach infinity, because $\frac{\partial L}{\partial
\beta}$ then approaches zero. This is of course impossible, since {\tmem{z}}
is confined to a (bounded) gap. Therefore one has to replace the real
parameter {\tmem{z}} by a complex one, $z = t + i \varepsilon$, for a small
positive number $\varepsilon$. Repeating the steps just outlined in this
particular situation, and then letting $\varepsilon$ approach zero, yields the
claimed formula for the partial derivative $\frac{\partial L}{\partial
\beta}$.

2) It is worthwhile mentioning, that at this stage it is already clear that
any possible critical point for the function {\tmem{L}} has to be a local
maximum. Indeed, taking the second partial derivatives
\[ \frac{\partial^2 L}{\partial z^2} (\beta, z) = - \tau ((z - h (\beta))^{-
   2}), \]
\[ \frac{\partial^2 L}{\partial \beta \partial z} (\beta, z) = \tau ((z - h
   (\beta))^{- 1} (v + v^{\ast}) (z - h (\beta))^{- 1}) = \tau ((z - h
   (\beta))^{- 2} (v + v^{\ast})), \]
\[ \frac{\partial^2 L}{\partial \beta^2} (\beta, z) = - \tau ([(z - h
   (\beta))^{- 1} (v + v^{\ast})]^2), \]
and applying the Cauchy-Schwarz inequality shows that the determinant of the
Hessian of the function {\tmem{L}} is strictly positive. Since the diagonal
entries of the Hessian are negative numbers, the claim follows.

.

{\tmstrong{Section 2}}

We now turn to the expansion and refinement of the settings introduced in
Section 1. In the sequel we assume throughout that $\alpha$ is a fixed
irrational number, and that $\beta$ is a positive number less than 1. For $p,
q \epsilon \mathbb{Z}$, we define the standardized monomials
\[ w_{p q} = e^{- p q \pi \alpha i} u^p v^q, \]
and for $z \epsilon \mathbb{C} \backslash S p (h (\beta))$,
\[ c_{p q} (z) = \tau ((h (\beta) - z)^{- 1}) w_{p q} . \]
The standardization ensures that these numbers are real valued whenever
{\tmem{z}} is a real number. The double sequence $\{ c_{p q} (z) \}$ solves
the following system linear equations for $s = z$,

(2.1)
\[ \cos \pi \alpha q (x_{p + 1, q} + x_{p - 1, q}) + \beta \cos \pi \alpha p
   (x_{p, q + 1} + x_{p, q - 1}) = s x_{p q}, \]
\[ \sin \pi \alpha q (x_{p + 1, q} - x_{p - 1, q}) - \beta \sin \pi \alpha p
   (x_{p, q + 1} - x_{p, q - 1}) = 0, \]
for all $p, q \epsilon \mathbb{Z}$, except $p = q = 0$. We shall refer to the
system (2.1) with the case \ $p = q = 0$ exempted by (2.1)\tmrsup{$\ast$}.

Before we proceed with our objective, we are going to dwell a little on the
linear system (2.1). First, if one multiplies the second equation by the
imaginary unit {\tmem{i}}, and adds the result to the first equation, then one
obtains Harper's equation on the two dimensional lattice, which is so common
in physics, and there are no algebraic complications attached to this
equation. By contrast, the system (2.1) and its truncated version
(2.1)\tmrsup{$\ast$}, combines two features which call for a specifically
designed approach to generate and analyze its solutions (see [R3]). On the one
hand, the system is largely overdetermined, giving rise to redundancies:
Asymptotically there are roughly twice as many equations than variables. On
the other hand, the system is degenerate along the diagonals $p = q$ and $p =
- q$: A recursion involving $4 \times 4$ matrices to generate the solutions of
this system collapses as one tries to cross either one of those two axes. The
tenet underlying the proof of Proposition 1.3 is that these two features,
reflecting intrinsic properties of the operator $h (\beta$), must also hold
the key to its more elusive spectral properties. Generally speaking, for
numbers {\tmem{s}} in the spectrum of $h (\beta$) the system (2.1) yields
uniformly bounded solutions, which are obtained by evaluating certain states
defined on the $C^{\ast}$-algebra $A_{\alpha}$ at the standardized monomials
$w_{p q}$. (In [R1], where it was shown that the linear dimension of the space
of uniformly bounded solutions is always either equal to one or to two, they
were referred to as ``eigenstates''.) For {\tmem{s=z}} in the resolvent set of
$h (\beta$), the double sequence $\{ c_{p q} (z) \}$ is always exponentially
decaying as $| p | \rightarrow \infty$ and \ $| q | \rightarrow \infty$. This
is the crucial property we shall exploit in the proof. However, in oder to
take advantage of it, we first need to ``homogenize'' the double sequence. In
other words, we have to find a double sequence solving the homogeneous system
(2.1), which preserves some measure of that exponential decay, but also shares
the vanishing conditions that hail from a possible critical point for the
function {\tmem{L}}. To this effect we need to delve a bit deeper into system
(2.1)\tmrsup{$\ast$}.

The solutions to this system form a linear space of dimension 6. Up to a
scaling factor to be determined below, there are four solutions with the
property that each of those has non-vanishing coefficients only in exactly one
of the four sectors separated by the lines $p = q$ and $p = - q$ in the two
dimensional lattice. The components of these four solutions are nothing but
the ``Fourier coefficients'' of the resolvent of perturbations of $h (\beta)$
which are obtained by multiplying the generators {\tmem{u}} and {\tmem{v}}
with suitable complex numbers of modulus larger than one or less than one. Put
in technical terms, expanding these resolvents in terms of the standardized
monomials $w_{p q}$ yields a ``non-commutative'' multiple Laurent series whose
coefficients are exactly those solutions. To assign a symbol to each of the
four solutions, let's say that $R^{(1, 0)}_{p q}$(s) vanishes for $p \leqslant
0, R^{(- 1.0)}_{p q} (s$) vanishes for $p \geq 0, R^{(0, 1)}_{p q} (s$)
vanishes for $q \leq 0$, and $R^{(0, - 1)}_{p q} (s$) vanishes for $q \geq 0$.
All four of these solutions can be computed by means of a two component
recursion of the type (3.13) in [R3]. Since we assumed $\beta$ to be less than
one, it follows that $R^{(1, 0)} $(s) and $R^{(- 1, 0)} $(s) decay
exponentially of order $\beta$ along any line in the lattice with slope 1 or
$- 1$, and that $R^{(0, 1)}  (s$) as well as $R^{(0, - 1)}  (s$) grow
exponentially of order $\beta^{- 1}$ along the four lines with slope 1 or $-
1$ through the points (0,1) \ or ($0, - 1$) in exactly one direction. \ Taking
the arithmetic mean of $R^{(1, 0)} $({\tmem{z}}) and $R^{(- 1,
0)}$({\tmem{z}}), both suitably scaled, yields a solution $\{ d_{p q} (z) \}$
of (2.1)\tmrsup{$\ast$} which has the following properties

(2.2)
\[ d_{p q} (z) = d_{| p | | q |} (z), \tmop{for} \tmop{all} p, q \epsilon
   \mathbb{Z}; d_{1, 0} (z) = \frac{1}{2} ; d_{p q} (z) = 0 \tmop{for}
   \begin{array}{l}
     | q | \geq | p |
   \end{array} . \]
(2.3)
\[ \overline{\lim_{| p | \rightarrow \infty}} \beta^{- | p |} \begin{array}{c}
     | d_{p, k + p} (z) |
   \end{array} < \infty, \overline{\lim_{| p | \rightarrow \infty}} \beta^{- |
   p |} \begin{array}{c}
     | d_{p, k - p} (z) |
   \end{array} < \infty, \tmop{for} \tmop{all} k \epsilon \mathbb{Z}. \]
Returning to our objective, it now follows that the double sequence $\phi_{p
q} (z) = c_{p q} (z) - d_{p q} (z)$ solves the system (2.1) and has the
following additional properties:

(2.4)
\[ \{ \phi_{p q} (z) \} \tmop{decays} \tmop{two} - \tmop{sided}
   \tmop{exponentially} \tmop{along} \tmop{any} \tmop{line} \tmop{in}
   \tmop{the} \tmop{lattice} \tmop{with} \tmop{slope} 1 \tmop{or} - 1 \]
(2.5)
\[ \tmop{If} c_{00} (z) = c_{01} (z) = 0, \tmop{then} \phi_{p q} (z) = 0
   \tmop{for} | p |, | q | \leq 1 \]
We are going to shelve this for a while, and turn to the discussion of a
certain subalgebra of the $C^{\ast}$-algebra $A_{\alpha}$. For two elements
$a, b \epsilon A_{\alpha}$ we denote by $a l g^{\star} (a, b)$ the $\ast
$-algebra generated by {\tmem{a}} and {\tmem{b}}. Next we define two
distinguished elements.
\begin{eqnarray*}
  \mathbb{U}= \beta^{- \frac{1}{2}} u + \beta^{\frac{1}{2}} v, \mathbb{V}=
  w_{- 1, 1} &  & 
\end{eqnarray*}
which satisfy the relations

(2.6)
\[ \mathbb{U}\mathbb{V}= \lambda^{- 2} \mathbb{V}\mathbb{U}, \mathbb{U}^{\ast}
   \mathbb{V}= \lambda^2 \mathbb{V}\mathbb{U}^{\ast}, \tmop{where} \lambda =
   e^{\pi \alpha i} . \]
One way to look at these two elements is that, while mimicking the generators
for the rotation algebra $A_{\alpha}$, they also allow for the representation
of the element $h (\beta)$ in a form that resembles representing the
degenerate element {\tmem{h}}(0) in terms of {\tmem{u}} and {\tmem{v}},
\[ h (\beta) =\mathbb{U}+\mathbb{U}^{\ast} . \]
The next step is to complete the mimicry, rendering $a l g^{\star}
(\mathbb{U}, \mathbb{V})$ as much as possible a look-alike of $a l g^{\star}
(u, v)$. What's missing from the picture is a basic symmetry. It is a basic
property of the irrational rotation algebra that there are two distinct
involutions that assign to one of the generators its adjoint, while fixing the
other one. One is a linear antiautomorphism, which we shall denote $\chi$, the
other one is a conjugate linear automorphism which we shall denote $\sigma$.

(2.7)
\[ \chi (u) = \sigma (u) = u^{\ast}, \chi (v) = \sigma (v) = v . \]
While both of these symmetries transform $a l g^{\star} (u, v)$ into itself,
this is not true for $a l g^{\star} (\mathbb{U}, \mathbb{V})$, due to the
asymmetric nature of the elements $\mathbb{U}$ and $\mathbb{V}$. To obtain two
such symmetries for $a l g^{\star} (\mathbb{U}, \mathbb{V})$, we first
introduce an automorphism of $A_{\alpha}$ that appeared for the first time in
[R4],

(2.8)
\[ \rho_{\beta} (u) = v u v (u v + \beta)^{- 1} (v^{\ast} u^{\ast} + \beta) \]
\[ \rho_{\beta} (v) = v (u v + \beta)^{- 1} (v^{\ast} u^{\ast} + \beta) . \]
Notice that, due to the general properties of $A_{\alpha}$, it is quite easy
to define automorphisms of $A_{\alpha}$. All one has to do is to assign
unitary elements in $A_{\alpha}$ to the two generators which preserve the
fundamental commutation relation (1.1) such that their images are also
generators. Any assignment of this kind extends automatically to an
automorphism. Since it can be shown that $\mathbb{U}$ is a generator for
$A_{\alpha}$, in other words the set of all (non-commutative) polynomials in
$\mathbb{U}$ and $\mathbb{U}^{\ast}$ is norm dense in $A_{\alpha}$, the
automorphism $\rho_{\beta}$ is uniquely determined by the identity,

(2.9)
\[ \rho_{\beta} (u + \beta v) = u^{\ast} + \beta v . \]
Thus, composition of $\chi$ and $\sigma$ \ with $\rho_{\beta}$,
\[ \chi_{\beta} = \chi  \circ \rho_{\beta}, \quad \sigma_{\beta} = \sigma
   \circ \rho_{\beta}, \]
yields the desired symmetries which are uniquely determined by the assignments

(2.10)
\[ u \rightarrow v u^{\ast} v (u^{\ast} v + \beta)^{- 1} (v^{\ast} u + \beta),
\]
\[ v \rightarrow v (u^{\ast} v + \beta)^{- 1} (v^{\ast} u + \beta) . \]
\

The restriction of these symmetries to the algebra $a l g^{\star}
(\mathbb{U}, \mathbb{V})$ share the desired properties,

(2.11)
\[ \chi_{\beta} (\mathbb{U}) =\mathbb{U}, \chi_{\beta} (\mathbb{V})
   =\mathbb{V}^{\star}, \qquad \sigma_{\beta} (\mathbb{U}) =\mathbb{U},
   \sigma_{\beta} (\mathbb{V}) =\mathbb{V}^{\star} . \]
While the first relation is obvious, the second one can be checked through
straightforward manipulations. Note that the symmetries $\chi_{\beta}$ and
$\sigma_{\beta}$ evaluated at standard monomials are indistinuishable, they
take the same values. Our next objective is to show that these symmetries,
when evaluated at the standardized monomials, yield elements whose expansion
in the standardized monomials have desirable exponential decay properties.
First we observe that
\[ (u^{\ast} v + \beta)^{- 1} = \sum_{n = 0}^{\infty} (- \beta)^n (v^{\ast}
   u)^{n + 1}, \]
and taking the adjoint on both sides yields of course a similar expansion.
This shows that each symmetry applied to $w_{p q}$ is a product of $v^{p + q}$
and an element that has a power series expansion in the monomial $w_{1, - 1}$
whose radius of convergence is equal to $\beta^{- 1}$. In conclusion, we
obtain the following representation,

(2.12)
\[ \sigma_{\beta} (w_{p q}) =  \sum_{m \epsilon \mathbb{Z}} r_m^{(p, q)} w_{m,
   p + q - m}, \tmop{where} \overline{\lim_{| m | \rightarrow \infty}}  
   \begin{array}{l}
     \begin{array}{l}
       |
     \end{array} r_m^{(p, q)} |^{\frac{1}{| m |}} \leqslant \beta
   \end{array} . \]
In preparation for the proof of Proposition 1.3 in the next section we
introduce two linear functional $\varphi_z$ \ and $\varphi_z \bullet
\sigma_{\beta}$ which are defined for $z \epsilon \mathbb{R} \backslash S p (h
(\beta))$ on the algebra $a l g^{\ast} (u, v)$ by the assignments
\[ \varphi_z (w_{p q}) = \phi_{p q} (z), \]
\[ \varphi_z \bullet \chi_{\beta} (w_{p q}) = \sum_{m \epsilon \mathbb{Z}}
   r_m^{(p, q)} \varphi_z (w_{m, p + q - m}), \quad \varphi_z \bullet
   \sigma_{\beta} (w_{p q}) = \sum_{m \epsilon \mathbb{Z}} \overline{r_m^{(p,
   q)} \varphi_z (w_{m, p + q - m})} \]
Notice that, by (2.4) and (2.12) the terms in the sum on the right-hand side
decay exponentially of an order less than or equal to $\beta$. These last
formulas define essentially the composition of the \ functional $\tmop{with}
\tmop{the} \tmop{two} \tmop{symmetries}$. Since $\sigma_{\beta}$ is conjugate
linear, we have to conjugate the terms in the sum in order to render the
resulting functional linear. We are now going to show that the two functionals
are actually equal.

(2.13)
\[ \varphi_z = \varphi_z \bullet \chi_{\beta}, \quad \varphi_z = \varphi_z
   \bullet \sigma_{\beta} . \]
Since the case of $\chi_{\beta}$ is slightly technically simpler because there
are no conjugates involved, we will show the arguments for $\sigma_{\beta}$
only. To see that the two functionals are the same, we need to return briefly
to the settings at the beginning of the present section. First, since
$\sigma_{\beta} (h (\beta)) = h (\beta)$, we also have $\sigma_{\beta} ((h
(\beta) - z)^{- 1}) = (h (\beta) - z)^{- 1}$. Notice that, by our assumption,
$z$ is a real number. This means that all we need to show is, that the
functionals $\vartheta_z$ and $\vartheta_z \bullet \sigma_{\beta}$, which are
defined below, are equal.

(2.14)
\[ \vartheta_z (w_{p q}) = R^{(1, 0)}_{p q} (z), \vartheta_z \bullet
   \sigma_{\beta} (w_{p q}) = \sum_{m \epsilon \mathbb{Z}} \overline{r_m^{(p,
   q)} \vartheta_z (w_{m, p + q - m})} . \]
and that a similar statement holds for $R^{(- 1, 0)}  (z)$. Again, since
$\sigma_{\beta} (h (\beta)) = h (\beta)$, and since $R^{(1, 0)}$({\tmem{z}}) \
solves the system (2.1)\tmrsup{$\ast$}, the second functional solves the
system (2.1)\tmrsup{$\ast$} as well. However, since by the representation in
(2.12) and the vanishing properties of $R^{(1, 0)}$({\tmem{z}}), $\vartheta_z
\bullet \sigma_{\beta} (w_{p q}$) vanishes for all \ indices located below or
on the the line $p = - q$ in the two dimensional lattice, the double sequence
$\{ \vartheta_z \bullet \sigma_{\beta} (w_{p q}) \}$ must be a linear
combination of $R^{(1, 0)}$({\tmem{z}}) and $R^{(0, 1)}$({\tmem{z}}). Since
$R^{(0, 1)}$({\tmem{z}}) grows exponentially of order $\beta^{- 1}$ along the
two lines with slope 1 or$- 1$ through the point ($0, 1$), while \
$\vartheta_z \bullet \sigma_{\beta} (w_{p q}$) and $R_{p q}^{(1,
0)}$({\tmem{z}}) vanish or decay exponentially along these two lines, it
follows that $\{ \vartheta_z \bullet \sigma_{\beta} (w_{p q}) \}$ is just a
scalar multiple of $R^{(1, 0)}$({\tmem{z}}). But since $\sigma_{\beta}$ is
unital, the two double sequences must actually be equal. A similar statement
can now be obtained for $R^{(- 1, 0)}  (z$) along a similar line of reasoning.
This concludes our argument establishing the validity of (2.13).

{\tmstrong{Remark:}} Restricting the solutions of the system
(2.1)\tmrsup{$\ast$}, viewed as linear functionals on the linear space $a l
g^{\star} (u, v$), to the subspace $a l g^{\star} (\mathbb{U}, \mathbb{V}$),
reduces the linear dimension from 6 to 5.

.

{\tmstrong{Section 3}}

We turn now to the proof of Proposition 1.3. Henceforth we shall simply write
$\varphi$ \ for the functional $\varphi_z$ introduced in section 2, because we
shall \ assume that $z$ is a fixed real number in the resolvent set of $h
(\beta)$. The idea of the proof is to exploit the decay conditions of $\varphi
(w_{p q})$ along lines with slope $- 1$in the two dimensional lattice, to
construct a functional $\psi$ \ on $a l g^{\star} (u, v)$, with the property
that $\varphi$ can be recovered from $\psi$ by the identity $\psi
(a\mathbb{U}) = \varphi (a)$ for all $a \epsilon a l g^{\star} (u, v)$, and
then showing that such a functional can not exist, in case the function
{\tmem{L}} in Section 1 has a critical point. Of course, there are infinitely
many functionals with this property. One simply has to implement a separate
elementary recursion along every single line with slope $- 1$. The crux is to
impose constraints which restrict the availability of such functionals
severely. More specifically, we shall prove the following.

\

{\tmstrong{3.1 Lemma}} There exists a linear functional $\psi$ on $a l
g^{\star} (u, v)$, having the properties,
\[ (i) \forall a \epsilon a l g^{\star} (u, v) : \psi (a\mathbb{U}) = \varphi
   (a) . \]
\[ (\tmop{ii}) \forall p \epsilon \mathbb{N}_0, \forall q \epsilon \mathbb{Z}:
   \psi (\mathbb{U}^p \mathbb{V}^q) = \overline{\psi (\mathbb{U}^p
   \mathbb{V}^{- q}}), \psi ((\mathbb{U}^{\ast})^p \mathbb{V}^q) =
   \overline{\psi ((\mathbb{U}^{\ast})^p \mathbb{V}^{- q}}) . \]
\[ (\tmop{iii}) \forall a \epsilon a l g^{\star} (u, v) : \varphi ((h (\beta)
   - z) a) = \psi (a (h (\beta) - z)\mathbb{U}) = 0 . \]

{\tmstrong{Proof:}} First we extend slightly the settings of Section 2. Let
$\mathfrak{B}$ be the set of elements {\tmem{a}} in $A_{\alpha}$ which can be
written in the form

(3.1)
\[ a = \sum^N_{n = - N} \sum_{m \epsilon \mathbb{Z}} k_{m, n - m} w_{m, n -
   m}, \tmop{where} \overline{\lim_{| m | \longrightarrow \infty}}
   \begin{array}{l}
     | k_{m, n - m} |
   \end{array}^{\frac{1}{| m |}} < \infty \tmop{for} - N \leq n \leq N . \]
This set is an involutive subalgebra of $A_{\alpha}$. Also, (2.12) implies

(3.2)
\[ \chi_{\beta} (\mathfrak{B}) =\mathfrak{B}, \quad \sigma_{\beta}
   (\mathfrak{B}) =\mathfrak{B}. \]
Furthermore,
\[ \mathbb{U}^{- 1} = \beta^{- \frac{1}{2}} \sum^{\infty}_{n = 0} (- \beta)^n
   (u^{\ast} v)^n u^{\ast} \epsilon \mathfrak{B}. \]

Now (2.3) allows us to extend the definition of the functional $\varphi$ to
elements of the form (3.1) as follows,
\[ \varphi (a) = \sum^N_{n = - N} \sum_{m \epsilon \mathbb{Z}} k_{m, n - m}
   \varphi (w_{m, n - m}) . \]

Obviously, this extended functional is also linear. Moreover, by (2.13)

(3.3)
\[ \varphi (\chi_{\beta} (a)) = \varphi (a), \quad \varphi (\sigma_{\beta}
   (a)) = \overline{\varphi (a)}, a \epsilon \mathfrak{B}, \]
and since the double sequence $\{ \phi_{p q} (z) \}$ in Section 2 solves the
system (2.1),

(3.4)
\[ \varphi ((h (\beta) - z) a) = \varphi (a (h (\beta) - z)) = 0, a \epsilon
   \mathfrak{B}. \]
We are now going to define the functional $\psi$ as follows,
\[ \psi (a) = \varphi (a\mathbb{U}^{- 1}), a \epsilon \mathfrak{B}. \]
We need to check that $\psi$ has the claimed properties. By construction this
is obvious for (i). Next, if {\tmem{p}} is a non-negative integer, and
{\tmem{q}} is an arbitrary integer, then (2.11) and (3.3) yield on the one
hand,
\[ \psi (\mathbb{U}^p \mathbb{V}^q) = \varphi (\mathbb{U}^p \mathbb{V}^q
   \mathbb{U}^{- 1}) = \overline{\varphi (\sigma_{\beta} (\mathbb{U}^p
   \mathbb{V}^q \mathbb{U}^{- 1}))} = \overline{\varphi (\mathbb{U}^p
   \mathbb{V}^{- q} \mathbb{U}^{- 1})} = \overline{\psi (\mathbb{U}^p
   \mathbb{V}^{- q})}, \]
and on the other hand,
\[ \begin{array}{ll}
     \psi (\mathbb{U}^p \mathbb{V}^q) & = \varphi (\mathbb{U}^p \mathbb{V}^q
     \mathbb{U}^{- 1}) = \lambda^{- 2 q} \varphi (\mathbb{U}^p \mathbb{U}^{-
     1} \mathbb{V}^q) = \lambda^{- 2 q} \varphi (\chi_{\beta} (\mathbb{U}^p
     \mathbb{U}^{- 1} \mathbb{V}^q) =
   \end{array} \]
\[ \lambda^{- 2 q} \varphi (\mathbb{V}^{- q} \mathbb{U}^p \mathbb{U}^{- 1})
   \begin{array}{l}
     = \lambda^{- 2 q} \psi (\mathbb{V}^{- q} \mathbb{U}^p) = \lambda^{- 2 (q
     + p q)} \psi (\mathbb{U}^p \mathbb{V}^{- q}) . \quad
   \end{array} \]
Combining these two lines of reasoning establishes the first identity in (ii).
The second identity can be shown in the same way. Finally, (3.4) yields,
\[ \psi ((h (\beta) - z) a) = \varphi ((h (\beta) - z) a\mathbb{U}^{- 1}) =
   0, a \epsilon \mathfrak{B}, \]
and also,
\[ \psi (a (h (\beta) - z)\mathbb{U}) = \varphi (a (h (\beta) -
   z)\mathbb{U}\mathbb{U}^{- 1}) = \varphi (a (h (\beta) - z)) = 0, \]
which establishes (iii) as well. $\blacktriangleleft$

\

{\tmstrong{Remark:}} Tracking the significance of the decay condition (2.4)
through the discussion so far, one observes that this condition is far
stronger than what is needed to make the arguments work. It would be enough to
assume that the double sequence in (2.4) does not increase exponentially of an
order larger than or equal to $\beta^{- 1}$ along any line with slope $- 1$ \
in the two dimensional lattice. All one has to do is to impose a stronger
exponential decay condition on the elements in the algebra $\mathfrak{B}$.

\

In the proof of the following lemma we shall be using nothing but the
relation (2.6), as well as the following

(3.5)
\[ \mathbb{U}^{\ast} \mathbb{U}= \lambda \mathbb{V}+ \lambda^{- 1}
   \mathbb{V}^{\ast} + \gamma \mathbb{e}, \tmop{where} \gamma = \beta +
   \beta^{- 1} \]

{\tmstrong{3.2 Lemma}} If $\psi$ is a linear functional defined on $a l
g^{\star} (\mathbb{U}, \mathbb{V})$ such that for some non-zero real number
{\tmem{t}} the following \ conditions hold,
\[ (i) \forall a \epsilon a l g^{\star} (\mathbb{U}, \mathbb{V}) : \psi
   ((\mathbb{U}+\mathbb{U}^{\ast} - t\mathbb{e}) a) = 0, \psi (a
   (\mathbb{U}+\mathbb{U}^{\ast} - t\mathbb{e})\mathbb{U}) = 0, \]
\[ (\tmop{ii}) \psi (\mathbb{U}) = \psi (\mathbb{U}^{\ast}) = \psi
   (\mathbb{V}) = \psi (\mathbb{V}^{\ast}) = \psi (\mathbb{e}) = 0, \]
\[ (\tmop{iii}) \forall p \in \mathbb{N}_0, \forall q \in \mathbb{Z}: \psi
   (\mathbb{U}^p \mathbb{V}^q) = \lambda^{- 2 (q + p q)} \overline{\psi
   (\mathbb{U}^p \mathbb{V}^q)}, \quad \forall p \in \mathbb{N}_0, \forall q
   \in \mathbb{Z}: \psi ((\mathbb{U}^{\star})^p \mathbb{V}^q) = \lambda^{- 2
   (q - p q)} \overline{\psi ((\mathbb{U}^{\star})^p \mathbb{V}^q),} \]
then $\psi (\mathbb{U}\mathbb{V}) = 0$.

\

{\tmstrong{Proof: }}First we derive a basic identity to be employed later in
the proof. Since
\[ (\mathbb{U}+\mathbb{U}^{\ast} - t\mathbb{e})\mathbb{U}\mathbb{V}^q
   =\mathbb{U}^2 \mathbb{V}^q +\mathbb{U}^{\ast} \mathbb{U}\mathbb{V}^q -
   t\mathbb{U}\mathbb{V}^q, \]
\[ \mathbb{V}^q (\mathbb{U}^2 +\mathbb{U}^{\ast} \mathbb{U}- t\mathbb{U}) =
   \lambda^{4 q} \mathbb{U}^2 \mathbb{V}^q +\mathbb{U}^{\ast}
   \mathbb{U}\mathbb{V}^q - \lambda^{2 q} t\mathbb{U}\mathbb{V}^q, \]
(i) yields
\[ \psi (\mathbb{U}^{\ast} \mathbb{U}\mathbb{V}^q) = t \psi
   (\mathbb{U}\mathbb{V}^q) - \psi (\mathbb{U}^2 \mathbb{V}^q), \]
\[ \psi (\mathbb{U}^{\ast} \mathbb{U}\mathbb{V}^q) = t \lambda^{2 q} \psi
   (\mathbb{U}\mathbb{V}^q) - \lambda^{4 q} \psi (\mathbb{U}^2 \mathbb{V}^q),
\]
hence,

(3.6)
\[ t \psi (\mathbb{U}\mathbb{V}^q) = (1 + \lambda^{2 q}) \psi (\mathbb{U}^2
   \mathbb{V}^q), \tmop{for} q \neq 0 . \]
Employing the second identity in (i) with $a =\mathbb{V}\mathbb{U}^{\ast}$
yields,
\[ \psi (\mathbb{V}\mathbb{U}^{\ast} (\mathbb{U}^2 +\mathbb{U}^{\ast}
   \mathbb{U}- t\mathbb{U})) = 0, \]
hence
\[ \psi (\mathbb{V}(\mathbb{U}^{\ast} \mathbb{U})\mathbb{U}) + \psi
   (\mathbb{V}\mathbb{U}^{\ast} (\mathbb{U}^{\ast} \mathbb{U})) - t \psi
   (\mathbb{V}(\mathbb{U}^{\ast} \mathbb{U})) = 0, \]
which by virtue of (3.5) yields,
\[ \psi (\mathbb{V}(\lambda \mathbb{V}+ \lambda^{- 1} \mathbb{V}^{\ast} +
   \gamma \mathbb{e})\mathbb{U}) + \psi (\mathbb{V}\mathbb{U}^{\ast} (\lambda
   \mathbb{V}+ \lambda^{- 1} \mathbb{V}^{\ast} + \gamma \mathbb{e})) - t \psi
   (\mathbb{V}(\lambda \mathbb{V}+ \lambda^{- 1} \mathbb{V}^{\ast} + \gamma
   \mathbb{e})) = 0, \]
hence by (ii),
\[ \psi (\mathbb{V}^2 \mathbb{U}) + \gamma \psi (\mathbb{V}\mathbb{U}) +
   \lambda \psi (\mathbb{V}\mathbb{U}^{\ast} \mathbb{V}) + \gamma \psi
   (\mathbb{V}\mathbb{U}^{\ast}) - t \lambda \psi (\mathbb{V}^2) = 0 \]
or by (2.6),

(3.7)
\[ \lambda^5 \psi (\mathbb{U}\mathbb{V}^2) + \gamma \lambda^2 \psi
   \mathbb{U}\mathbb{V}) + \lambda^{- 1} \psi (\mathbb{U}^{\ast} \mathbb{V}^2)
   + \gamma \lambda^{- 2} \psi (\mathbb{U}^{\ast} \mathbb{V} ) - t \lambda
   \psi (\mathbb{V}^2) = 0. \]
Conjugating this equation, while exploiting the fact that the second identity
in (iii) ensures that $\psi (\mathbb{U}^{\ast} \mathbb{V} )$ and $\psi
(\mathbb{U}^{\ast} \mathbb{V}^2)$ are real numbers implies,

(3.8)
\[ \lambda^{- 5} \overline{\psi (\mathbb{U}\mathbb{V}^2)} + \gamma \lambda^{-
   2} \overline{\psi \mathbb{U}\mathbb{V})} + \lambda  \psi (\mathbb{U}^{\ast}
   \mathbb{V}^2) + \gamma \lambda^2 \psi (\mathbb{U}^{\ast} \mathbb{V} ) - t
   \lambda^{- 1} \overline{\psi (\mathbb{V}^2)} = 0 \quad \]
(3.7) and (3.8) form an inhomogeneous linear system in the two variables $\psi
(\mathbb{U}^{\ast} \mathbb{V} )$ and $\psi (\mathbb{U}^{\ast} \mathbb{V}^2)$:
\[ \left(\begin{array}{c}
     \begin{array}{l}
       \lambda^{- 1}\\
       \lambda
     \end{array} \begin{array}{l}
       \gamma \lambda^{- 2}\\
       \gamma \lambda^2
     \end{array}
   \end{array}\right) \left( \begin{array}{l}
     \psi (\mathbb{U}^{\ast} \mathbb{V}^2)\\
     \psi (\mathbb{U}^{\ast} \mathbb{V} )
   \end{array} \right) = \left( \begin{array}{l}
     t \lambda \psi (\mathbb{V}^2) - \lambda^5 \psi (\mathbb{U}\mathbb{V}^2) -
     \gamma \lambda^2 \psi (\nobracket \mathbb{U}\mathbb{V})\\
     t \lambda^{- 1} \overline{\psi (\mathbb{V}^2)} - \lambda^{- 5}
     \overline{\psi (\mathbb{U}\mathbb{V}^2)} - \gamma \lambda^{- 2}
     \overline{\psi (\nobracket \mathbb{U}\mathbb{V})}
   \end{array} \right) . \]

Applying the first identity in (iii) to the terms in the second component on
the right hand side yields,
\[ \left(\begin{array}{c}
     \begin{array}{l}
       \lambda^{- 1}\\
       \lambda
     \end{array} \begin{array}{l}
       \gamma \lambda^{- 2}\\
       \gamma \lambda^2
     \end{array}
   \end{array}\right) \left( \begin{array}{l}
     \psi (\mathbb{U}^{\ast} \mathbb{V}^2)\\
     \psi (\mathbb{U}^{\ast} \mathbb{V} )
   \end{array} \right) = \left( \begin{array}{l}
     t \lambda \psi (\mathbb{V}^2) - \lambda^5 \psi (\mathbb{U}\mathbb{V}^2) -
     \gamma \lambda^2 \psi (\nobracket \mathbb{U}\mathbb{V})\\
     t \lambda^3 \psi (\mathbb{V}^2) - \lambda^3 \psi (\mathbb{U}\mathbb{V}^2)
     - \gamma \lambda^2 \psi (\nobracket \mathbb{U}\mathbb{V})
   \end{array} \right), \]
which in turn leads to the solution,

(3.9)
\[ \left( \begin{array}{l}
     \psi (\mathbb{U}^{\ast} \mathbb{V}^2)\\
     \psi (\mathbb{U}^{\ast} \mathbb{V} )
   \end{array} \right) = \left( \begin{array}{l}
     - \lambda (- \lambda t \psi (\mathbb{V}^2) + (\lambda^5 + \lambda^3 +
     \lambda) \psi (\mathbb{U}\mathbb{V}^2) + (\gamma \lambda^2 + \gamma) \psi
     (\nobracket \mathbb{U}\mathbb{V}))\\
     \frac{\lambda^2  ((\lambda^3 + \lambda) \psi (\mathbb{U}\mathbb{V}^2) +
     \gamma \psi (\nobracket \mathbb{U}\mathbb{V}))}{\gamma}
   \end{array} \right) \quad \]

Employing the first identity in (i) with $a = \mathbb{V}$ yields,
\[ \psi ((\mathbb{U}+\mathbb{U}^{\ast} - t\mathbb{e})\mathbb{V}) = 0 \]
Hence by (ii),
\[ \psi (\mathbb{U}\mathbb{V}) = - \psi ((\mathbb{U}^{\ast} \mathbb{V}) . \]
Combining this with the second component in (3.9) yields,
\[ \lambda^2  (\lambda^3 + \lambda) \psi (\mathbb{U}\mathbb{V}^2) + \gamma
   (\lambda^2 + 1) \psi (\mathbb{U}\mathbb{V}) = 0, \]
which simplyfies to,

(3.10)
\[ \lambda^3 \psi (\mathbb{U}\mathbb{V}^2) + \gamma \psi
   (\mathbb{U}\mathbb{V}) = 0. \]
Next, starting over again, using the first identity in (i) with $a
=\mathbb{U}^2 \mathbb{V}$,
\[ \psi ((\mathbb{U}+\mathbb{U}^{\ast} - t\mathbb{e})\mathbb{U}^2 \mathbb{V})
   = 0, \]
or equivalently,
\[ \psi (\mathbb{U}^3 \mathbb{V}) + \psi ((\mathbb{U}^{\ast} \mathbb{U}
   )\mathbb{U}\mathbb{V}) - t \psi (\mathbb{U}^2 \mathbb{V}) = 0, \]
yielding together with (3.5),
\[ \psi (\mathbb{U}^3 \mathbb{V}) + \psi ((\lambda \mathbb{V}+ \lambda^{- 1}
   \mathbb{V}^{\ast} + \gamma \mathbb{e})\mathbb{U}\mathbb{V}) - t \psi
   (\mathbb{U}^2 \mathbb{V}) = 0, \]
which by (ii) implies,
\[ \psi (\mathbb{U}^3 \mathbb{V}) + \lambda^3 \psi (\mathbb{U}\mathbb{V}^2) +
   \gamma \psi (\mathbb{U}\mathbb{V}) - t \psi (\mathbb{U}^2 \mathbb{V}) = 0 .
\]
Employing (3.6) to this equation for $q = 1$, we obtain,

(3.11)
\[ \psi (\mathbb{U}^3 \mathbb{V}) + \lambda^3 \psi (\mathbb{U}\mathbb{V}^2) +
   (\gamma - \frac{t^2}{1 + \lambda^2}) \psi (\mathbb{U}\mathbb{V}) = 0 . \]
On the other hand, the second identity in (i) with $a =\mathbb{U}\mathbb{V}$
yields,
\[ \psi (\mathbb{U}\mathbb{V}(\mathbb{U}^2 +\mathbb{U}^{\ast} \mathbb{U}-
   t\mathbb{U})) = 0, \]
or equivalently
\[ \psi (\mathbb{U}\mathbb{V}\mathbb{U}^2) + \psi
   (\mathbb{U}\mathbb{V}(\mathbb{U}^{\ast} \mathbb{U})) - t \psi
   (\mathbb{U}\mathbb{V}\mathbb{U}) = 0, \]
hence by (3.5),
\[ \psi (\mathbb{U}\mathbb{V}\mathbb{U}^2) + \psi
   (\mathbb{U}\mathbb{V}(\lambda \mathbb{V}+ \lambda^{- 1} \mathbb{V}^{\ast} +
   \gamma \mathbb{e})) - t \psi (\mathbb{U}\mathbb{V}\mathbb{U}) = 0, \]
which by (ii) and (2.6) yields,
\[ \lambda^4 \psi (\mathbb{U}^3 \mathbb{V}) + \lambda \psi
   (\mathbb{U}\mathbb{V}^2) + \gamma \psi (\mathbb{U}\mathbb{V}) - t \lambda^2
   \psi (\mathbb{U}^2 \mathbb{V}) = 0, \]
and after invoking (3.6) for $q = 1$ again,

(3.12)
\[ \lambda^4 \psi (\mathbb{U}^3 \mathbb{V}) + \lambda \psi
   (\mathbb{U}\mathbb{V}^2) + (\gamma - \frac{\lambda^2 t^2}{1 + \lambda^2})
   \psi (\mathbb{U}\mathbb{V}) = 0 . \]
(3.11) and (3.12) form a linear sytem in the three variables $\psi
(\mathbb{U}\mathbb{V})$, $\psi (\mathbb{U}\mathbb{V}^2)$, $\psi
(\mathbb{U}\mathbb{V}^2)$ with a one-dimensional null space generated by the
vector,
\[ \left( \frac{t^2 - \gamma}{(\lambda^2 - \lambda + 1)  (\lambda^2 + \lambda
   + 1)}, - \frac{\gamma \lambda^4 + 2 \gamma \lambda^2 + \gamma - \lambda^2
   t^2}{\lambda (\lambda^2 + 1)  (\lambda^2 - \lambda + 1)  (\lambda^2 +
   \lambda + 1)}, 1 \right) \]
Scaling and rearranging the component terms yields a vector with real entries
in the first and second position:
\[ \left( t^2 - \gamma, \frac{t^2}{\lambda + \bar{\lambda}} - (\lambda +
   \bar{\lambda}) \gamma, (\lambda^2 - \lambda + 1)  (\lambda^2 + \lambda + 1)
   \right) \]
There exists a complex number $z$ such that,

(3.13)
\[ z (t^2 - \gamma) = \psi (\mathbb{U}^3 \mathbb{V}), \quad z \left(
   \frac{t^2}{\lambda + \bar{\lambda}} - (\lambda + \bar{\lambda}) \gamma
   \right) = \psi (\mathbb{U}\mathbb{V}^2), \quad z (\lambda^2 - \lambda + 1) 
   (\lambda^2 + \lambda + 1) = \psi (\mathbb{U}\mathbb{V}) . \]

Combining (3.10) with (3.13) we conclude,
\[ \lambda^3 \psi (\mathbb{U}\mathbb{V}^2) + \gamma \psi
   (\mathbb{U}\mathbb{V}) = \lambda^3 z \left( \frac{t^2}{\lambda +
   \bar{\lambda}} - (\lambda + \bar{\lambda}) \gamma \right) + \gamma z
   (\lambda^2 - \lambda + 1)  (\lambda^2 + \lambda + 1) = 0 \]
If \ $\psi (\mathbb{U}\mathbb{V}) \neq 0$, then \ $z \neq 0$ as well. Hence,
\[ \lambda^3 \left( \frac{t^2}{\lambda + \bar{\lambda}} - (\lambda +
   \bar{\lambda}) \gamma \right) + \gamma (\lambda^2 - \lambda + 1) 
   (\lambda^2 + \lambda + 1) = \frac{\gamma \lambda^2 + \gamma + \lambda^4
   t^2}{\lambda^2 + 1} = 0, \]
or equivalently,
\[ \frac{t^2}{\gamma} = - \frac{\lambda^2 + 1}{\lambda^4} . \]
Since $\frac{t^2}{\gamma}$ is a real number,
\[ - \frac{\lambda^2 + 1}{\lambda^4} - \left( - \frac{\lambda^{- 2} +
   1}{\lambda^{- 4}} \right) = 0, \]
or equivalently,
\[ \frac{(\lambda - 1) (\lambda + 1) (\lambda^2 + 1)  (\lambda^2 - \lambda +
   1)  (\lambda^2 + \lambda + 1)}{\lambda^4} = 0. \]
Since $\lambda$is not a root of unity, we have reached a contradiction.
Therefore, $\psi (\mathbb{U}\mathbb{V}) = 0$, as claimed.

\

{\tmstrong{Remark:}} Performing the kind of manipulations in the proof of
Lemma 2.3 for more general terms of the form $\mathbb{U}^p \mathbb{V}^q$ and
${(\mathbb{U}^{\ast}}^p)\mathbb{V}^q$ one can actually show that the
functional $\psi$ ``almost'' vanishes. If the condition $\psi
(\mathbb{U}^2)$=0 is added, then $\psi$ vanishes completely.

\

{\tmstrong{Proof of Proposition 1.3:}} First notice that the functional
$\varphi$ can not be zero. This is true because the double sequence $\{ c_{p
q} (z) \}$ decays exponentially as $| p | \rightarrow \infty$ and \ $| q |
\rightarrow \infty$, while this is not true for $\{ d_{p q} (z) \}$. If $\{
d_{p q} (z) \}$ were \ decaying exponentially as $| p | \rightarrow \infty$
and \ $| q | \rightarrow \infty$, then the two solutions $R_{1,
0}$({\tmem{z}}) and $R_{- 1, 0}$({\tmem{z}}) of the system
(2.1)\tmrsup{$\ast$} would give rise to two distinct inverses of the same
element $h (\beta) - z$, which is of course impossible. Now suppose that
\[ (\tau ((h (\beta) - z)^{- 1}) = \tau ((h (\beta) - z)^{- 1} v)) = 0 . \]
It follows from (2.5),
\[ \varphi (u) = \varphi (u^{\ast}) = \varphi (v) = \varphi (v^{\ast}) = 0 \]
Since the double sequence $\{ \varphi (w_{p q}) \}$ solves the system (2.1),
this implies that $\varphi$ is non-zero if and only if
\[ \varphi (\mathbb{V}) \neq 0 . \]
The vanishing conditions for $\varphi$ translate into several vanishing
conditions for $\psi$. First, by Lemma 3.1,
\[ \psi (\mathbb{U}) = \varphi (\mathbb{e}) = 0 . \]
Furthermore,
\[ \beta^{- \frac{1}{2}} \psi (\mathbb{e}) + \beta^{\frac{1}{2}} \psi
   (u^{\ast} v) = \psi (u^{\ast} \mathbb{U}) = \varphi (u^{\ast}) = 0, \]
\[ \beta^{- \frac{1}{2}} \psi (v^{\ast} u) + \beta^{\frac{1}{2}} \psi
   (\mathbb{e}) = \psi (v^{\ast} \mathbb{U}) = \varphi (v^{\ast}) = 0, \]
or, equivalently,
\[ \beta^{- \frac{1}{2}} \psi (\mathbb{e}) + \beta^{\frac{1}{2}} \lambda^{- 1}
   \psi (\mathbb{V}^{\ast}) = 0, \]
\[ \beta^{- \frac{1}{2}} \lambda \psi (\mathbb{V}) + \beta^{\frac{1}{2}} \psi
   (\mathbb{e}) = 0 . \]
By Lemma 3.1(ii),
\[ \psi (\mathbb{e}) = \overline{\psi (\mathbb{e})}, \psi (\mathbb{V}^{\ast})
   = \overline{\psi (\mathbb{V})} . \]
This yields a linear system for $\psi (\mathbb{e})$ and $\psi (\mathbb{V})$,
\[ \beta^{- \frac{1}{2}} \psi (\mathbb{e}) + \beta^{\frac{1}{2}} \lambda  \psi
   (\mathbb{V} ) = 0, \]
\[ \beta^{\frac{1}{2}} \psi (\mathbb{e}) + \beta^{- \frac{1}{2}} \lambda \psi
   (\mathbb{V}) = 0 . \]
Since the determinant of this system,
\[ \left|\begin{array}{c}
     \begin{array}{c}
       \beta^{- \frac{1}{2}}\\
       \beta^{\frac{1}{2}}
     \end{array} \begin{array}{c}
       \beta^{\frac{1}{2}} \lambda\\
       \beta^{- \frac{1}{2}} \lambda
     \end{array}
   \end{array}\right| = \lambda (\beta^{- 1} - \beta) \]
is non-zero for $\beta \neq 1$, it follows that
\[ \psi (\mathbb{e}) = \psi (\mathbb{V}) = \psi (\mathbb{V}^{\ast}) = 0 . \]
Finally, by Lemma 3.1(iii), since $\psi (\mathbb{e}) = 0$, and $\psi
(\mathbb{U}) = \varphi (\mathbb{e}) = 0$
\[ \psi (\mathbb{U}^{\ast}) = \psi (\mathbb{U}+\mathbb{U}^{\ast} -
   z\mathbb{e}) = 0 . \]
In conclusion we have shown, that $\psi$ (more precisely its restriction to $a
l g^{\star} (\mathbb{U}, \mathbb{V}$)) satisfies the conditions (i) and (ii)
of Lemma 2.3, letting $t = z \beta^{- \frac{1}{2}}$. Therefore,
\[ \varphi (\mathbb{V}) = \psi (\mathbb{V}\mathbb{U}) = \lambda^2 \psi
   (\mathbb{U}\mathbb{V}) = 0 . \]
Since we observed at the beginning of the proof that $\varphi (\mathbb{V})$
cannot be zero, we have reached a contradiction. $\blacktriangleleft$

.

{\tmstrong{References:}}

[BJ] J. Bourgain, S. Jitomirskaya, ``Continuity of the Lyapunov Exponent for
Quasiperiodic operators with analytic potential'', J. Stat. Phys. 108 (2002)
1203-1218

[E] H. M. Edwards, ``Riemann's Zeta Function'', Academic Press, 1974

[GSUKNKS] C. Geisler, J.H. Smet, V. Umansky, K. von Klitzing, B. Naundorf. R.
Ketzmerick H. Schweizer, ``Detection of a Landau Band-Coupling-Induced
Rearrangement of the \ Hofstadter Butterfly'', Phys. Rev. Lett. 92, 256801
(2004)

[HKS] B. Helffer, P.Kerdelhue, S. Sjostrand, ``La Papillon De Hofstadter
Revisite'', Supplement au Bulletin de la Societe Mathematique du France, 118,
3 (1990)

[JZ] D.Jaksch, P. Zoller, ``Creation of effective magnetic fields in optical
lattices: The Hofstadter butterfly for cold neutral atoms'', New J. Phys. 5
(2003) 56

[L] E. Landau, ``Vorlesungen Ueber Zahlentheorie'' Verlag von S. Hirzel,
Leipzig (1927)

[OA] D. Osadchy and J. E. Avron, ``Hofstadter butterfly as quantum phase
diagram'' J. Math. Phys. 42 (2001) 5665-5671

[R1] N. Riedel, ``Almost Mathieu operators and rotation $C^{\ast}
-$algebras'', Proc. London Math. Soc. 3 (56), (1988) 281-302

[R2] N. Riedel ``Regularity of the spectrum for the almost Mathieu operator'',
Proc. Amer. Math Soc., 129 (1999) 1681-1687

[R3] N. Riedel, ``The spectrum of a class of almost periodic operators'', Int.
J. Math. Sci. 36 (2003) 2277-2301

[R4] N.Riedel ``Exponentially decaying eigenvectors for certain almost
periodic operators'', Erg. Th. Dyn. Syst. 24 (2004) 915-943

.

E-mail: {\tmem{nriedel@tulane.edu}}

Department of Mathematics

Tulane University

New Orleans, LA 70118

\

\

\

\

\

\

\

\

\

\

\

\ \

\

\

\

\

\

\

\

\

\

\

\

\

\

\

\

\ \ \ \ \ \ \ \ \ \ \ \ \ \ \ \ \ \ \ \ \ \ \ \ \ \ \ \ \ \ \ \ \ \ \ \ \ \
\ \ \ \ \ \ \ \ \ \ \

\end{document}